\documentclass[11pt]{amsart}
\usepackage{hyperref}
\usepackage{a4wide,enumerate}
\usepackage{amsthm,amssymb}
\usepackage[english]{babel}
\newcommand{\RR}{\mathbb{R}}
\newcommand{\CC}{\mathbb{C}}
\newcommand{\NN}{\mathbb{N}}
\newcommand{\ZZ}{\mathbb{Z}}
\DeclareMathOperator{\EE}{\mathbb{E}}
\newcommand{\PP}{\mathbb{P}}

\newcommand{\cV}{\mathcal{V}}
\newcommand{\Tr}{{\mathop{\mathrm{Tr}}}}
\DeclareMathOperator{\vol}{vol}

\DeclareMathOperator{\supp}{\mathrm{supp}}

\newcommand{\dist}{{\ensuremath{\mathrm{dist}}}}

\DeclareMathOperator{\interior}{int}
\newcommand{\inter}[2]{\interior_{#2}{#1}}
\DeclareMathOperator{\exterior}{ext}
\newcommand{\exter}[2]{\exterior_{#2}{#1}}
\newcommand{\bdry}[2]{\partial_{#2}{#1}}
\def\incident{\rightsquigarrow}

\newcommand{\be}{\begin{equation}}
\newcommand{\ee}{\end{equation}}
\newcommand{\bea}{\begin{eqnarray*}}
\newcommand{\eea}{\end{eqnarray*}}

\newcommand{\Nach}{{\,\rightarrow\,}}

\newtheorem{thm}{Theorem}%[section]
\newtheorem{lem}[thm]{Lemma}

\newtheorem{cor}[thm]{Corollary}
\theoremstyle{definition}
\newtheorem{dfn}[thm]{Definition}

\theoremstyle{remark}
\newtheorem{rem}[thm]{Remark}
\newtheorem{exm}[thm]{Example}

\newcommand{\Hm}[1]{\leavevmode{\marginpar{\tiny%
$\hbox to 0mm{\hspace*{-0.5mm}$\leftarrow$\hss}%
\vcenter{\vrule depth 0.1mm height 0.1mm width \the\marginparwidth}%
\hbox to
0mm{\hss$\rightarrow$\hspace*{-0.5mm}}$\\\relax\raggedright #1}}}

\begin{document}
\title[Optimal Wegner estimates on metric graphs]
{Optimal Wegner estimates for random Schr\"odinger operators on metric graphs}

\author[M.~J.~Gruber]{Michael J.\ Gruber}
\address[M.G.]{TU Clausthal\\
Institut f\"ur Mathematik\\
38678 Clausthal-Zellerfeld\\
Germany}
\urladdr{\url{http://www.math.tu-clausthal.de/~mjg/}}

\author[M.~Helm]{Mario Helm}
\address[M.H.]{TU Chemnitz\\
Fakult\"at f\"ur Mathematik\\
09107 Chemnitz\\
Germany}

\author[I.~Veseli\'c]{Ivan Veseli\'c}
\address[I.V.]{Emmy-Noether-Programm der Deutschen Forschungsgemeinschaft \&
Fakul\-t\"at f\"ur Mathematik, 09107\, TU\, Chemnitz, Germany}
\urladdr{\url{http://www.tu-chemnitz.de/mathematik/enp}}
\curraddr[I.V.]{Institut f\"ur Angewandte Mathematik\\ 53115 Universit\"at Bonn\\ Germany}

\thanks{\copyright 2007 by the authors. Faithful reproduction of this article,
         in its entirety is permitted for non-commercial purposes. {\today, \jobname.tex}}



\begin{abstract}
We consider Schr\"odinger operators with a random potential of alloy type
on infinite metric graphs which obey certain uniformity conditions.
For single site potentials of fixed sign we prove that the random Schr\"odinger operator restricted to a
finite volume subgraph obeys a Wegner estimate which is linear in the volume and reproduces the
modulus of continuity of the single site distribution. This improves and unifies
earlier results for alloy type models on metric graphs.

We discuss applications of Wegner estimates to bounds on the modulus of continuity for
the integrated density of states of ergodic Schr\"odinger operators, as well
as to the proof of Anderson localisation via the multiscale analysis.
\end{abstract}

\maketitle
\let\languagename\relax

\bigskip
\section{Introduction}

Random Schr\"odinger operators have been extensively studied both on Euclidean space
$\RR^\nu$ and on the lattice $\ZZ^\nu$.
For a wide range of ergodic models several basic spectral properties
have been established.  These include the non-randomness of the spectrum
and its measure-theoretic components,
the existence of a self-averaging integrated density of states (IDS for short),
as well as a closed trace-per-unit-volume formula for the latter quantity.
Under more specific assumptions it has been proved that the IDS obeys certain continuity properties
and/or that the spectrum of the random operator is purely localised, at least in certain energy regions.
More precisely, for various models it has been shown that there exist energy intervals near
spectral boundaries and certain disorder regimes such that the random Schr\"odinger operator
exhibits pure point spectrum (in the mentioned energy interval) and that the
associated eigenfunctions decay exponentially, almost surely. This phenomenon is called spectral or \emph{Anderson localisation},
a term coined after the groundbreaking paper \cite{Anderson-58}.
Actually, even a stronger form of localisation holds for these operators, which is formulated in terms of the dynamical
properties of the time-evolution operator associated to the Schr\"odinger operator.
Since the literature on the mentioned models and results is vast, we refer only to the monographs
\cite{CarmonaL-90,PasturF-92,Stollmann-01} and the references therein.

More recently the spectral properties of similar models
on quantum respectively metric graphs have been analysed e.g.~in
\cite{KostrykinS-04,AizenmanSW-06a,HislopP,HelmV,GruberLV,ExnerHS,GruberV,KloppP}.

The present paper is devoted to the proof of a so called \emph{Wegner estimate} for rather general
random Schr\"odinger operators with non-negative alloy type potentials.
This type of estimate goes back to the paper \cite{Wegner-81} and
concerns the expected number of eigenvalues of the Schr\"odinger operator restricted to a finite volume
in a given energy interval. Let us stress certain interesting features of the models to which
our main result applies:
(1) \ The single-site potential needs to be positive on an open set,
but this set may be arbitrarily small.
(2) \ We do not need to assume a periodicity condition, since we do not rely on Floquet-Bloch theory.
(3) \ Our Wegner estimate reproduces the (arbitrary) modulus of continuity of the single site distribution.
In particular our results unify and extend the results of \cite{HelmV} and \cite{GruberV}.

The structure of the paper is as follows: In the second section we
introduce our model and state the Wegner estimate as the main
result. The third section describes two important applications of
such estimates, namely consequences for the modulus of continuity of
the IDS and localisation proofs via multiscale analysis, with an application
to $\log$-H\"older continuous single site distributions. The last
section contains proofs of our new key lemma and the main theorem, 
as well as a few lemmata taken from \cite{KirschV-02b},\cite{HelmV},
\cite{GruberV}, and \cite{GruberLV}, adapted to the new context:
arbitrary modulus of continuity, partial covering conditions.

M.G.\ acknowledges enlightening discussions with J.\ Brasche to whom we owe parts of Remark~\ref{r:rcmodulus}.
M.H.\ acknowledges partial support through the Deutsche Forschungsgemeinschaft.
I.V.\ acknowledges support through the Emmy-Noether programme of the Deutsche Forschungsgemeinschaft.

\section{Model and results}
\label{s-Results}

We start with the definition of metric graphs which are the
underlying topological structures of our models. A metric graph is a triple
$G=(V,E,\mathcal{G})$ of two countable sets $V$ and $E$ (vertices
and edges) and a map
\begin{equation*}
   \mathcal{G}: E\Nach V \times V \times (0,\infty), \qquad e\mapsto
   (\iota(e), \tau(e), l_e),
\end{equation*}
which determines for each edge $e$ an initial and terminal vertex and a
positive length $l_e$. In this sense we identify each edge $e$ with
the interval $(0, l_e)$. The pair $(V,E)$ is the combinatorial graph associated to
the metric graph $G=(V,E,\mathcal{G})$.  If all vertex degrees are finite we get a
metric on the topological space (CW-complex) $G$ by taking the infimum of
lengths of paths connecting two given points (see e.g.
\cite{Schubert-06}).

Given the metric graph $G$  we introduce induced
subgraphs. Here we write $v\incident e$ if $v\in V$ is incident to
$e\in E$.
\begin{dfn}\label{dfn-subgraph}
Let $(V,E)$ be a combinatorial graph without isolated vertices and $\Lambda\subset E$. Then we define a partition of $V$ as follows:
\begin{align*}
\inter V \Lambda &:= \big\{ v\in V : \{e\in E: v\incident e\}\subset \Lambda \big\}\quad\text{(interior vertices)} \\
\exter V \Lambda &:= \big\{ v\in V : \{e\in E: v\incident e\}\subset E\setminus\Lambda  \big\}\quad\text{(exterior vertices)} \\
\bdry  V \Lambda &:= \big\{ v\in V : \exists e\in\Lambda,e'\in  E\setminus\Lambda :  v\incident e,  v\incident e'  \big\}\quad\text{(boundary vertices)}
\end{align*}
With $V_\Lambda:=\inter V \Lambda \cup \bdry  V \Lambda $, the combinatorial subgraph induced by $\Lambda$ is then given  by $(V_\Lambda,\Lambda)$.
\end{dfn}
Note that $\bdry V {E\setminus\Lambda} = \bdry V \Lambda$ and $\inter V {E\setminus\Lambda} = \exter V \Lambda$.
The above partition of $V$ can equivalently be described by comparing vertex degrees, if we declare $\deg_{G_\Lambda}(v) = 0$ for $v\in V\setminus V_\Lambda$:
\begin{align*}
\inter V \Lambda &= \big\{ v\in V : \deg_{G_\Lambda}(v) = \deg_G(v) \big\} \\
\exter V \Lambda &= \big\{ v\in V : \deg_{G_\Lambda}(v) =0  \big\} \\
\bdry  V \Lambda &= \big\{ v\in V : 0 < \deg_{G_\Lambda}(v) <   \deg_G(v) \big\}
\end{align*}
Now we define the restriction of the metric graph $G=(V,E,\mathcal{G})$ induced by the edge-subset $\Lambda$ by
setting $G_\Lambda:=(V_\Lambda,\Lambda, \mathcal{G}|_{\Lambda})$. Here $\mathcal{G}|_{\Lambda}$
denotes the restriction of the map $\mathcal{G}$ to $\Lambda$.

Next we define the negative Laplacian on $L^2(G):=
\oplus_{e\in E} L^2(0,l_e)$. To this end we introduce the
Sobolev space $W_2^2(E):=\oplus_{e\in
E}W_2^2(0,l_e)$. Note that $W_2^2(0,l_e)\subset C^1([0,l_e])$, so that
the boundary values $f_e(0),f_e(l_e),f_e'(0),f_e'(l_e)$ are well defined
for $f\in W_2^2(E)$, where we set $f_e:=f|_e$.
We fix an arbitrary ordering for the edges $e$ incident to $v$
and write $f(v)\in\CC^{\deg v}$ for the vector of boundary values of $f_e$ at
$0$ (resp.\ $l_e$) if $\iota(e)=v$ (resp.\ $\tau(e)=v$).
Similarly, we write $\partial f(v)\in\CC^{\deg v}$ for the vector of
boundary values of $f_e'$ (resp.\ $-f_e'$) at
$0$ (resp.\ $l_e$) if $\iota(e)=v$ (resp.\ $\tau(e)=v$).

For each $v\in V$, let $S_v$ be a Lagrangian subspace of $\CC^{2\deg v}$ with respect to the standard complex symplectic structure
(see, e.g., \cite{KostrykinS-99b}).
Then the operator $-\Delta_G$ with boundary condition $(S_v)_{v\in V}$
is given by
\begin{eqnarray*}
   \mathcal D(-\Delta_G)&:=&\big\{ \,f\in W_2^2(E) \mid
    \forall v\in V: \big(f(v),\partial f(v)\big) \in S_v \big\},\\
   (-\Delta_G f)_e&:=& -f''_e \qquad(e\in E).
\end{eqnarray*}
This is the most general type of graph-local selfadjoint boundary conditions
and includes Dirichlet (Dirichlet on the edge), Neumann (Neumann on the edge),
% :)
free (Kirchhoff, Neumann, standard) boundary conditions, of course. We suppress the boundary conditions
in the notation since they will be fixed and clear from the context.

Note that all of the above applies to subgraphs as well, once we specify the
family of boundary value subspaces on the subgraph. These will be arbitrary
in our Wegner estimates, with constants independent of this choice!

On several occasions, we will make use of the so called \emph{restriction to $G_\Lambda$ with Dirichlet conditions}:
Given a Laplacian, i.e.\ a choice of boundary value spaces $S_v,v\in V,$ on $G$ and a subset $\Lambda\subset E$ we define
$\tilde S_v:=S_v$ for $v\in\inter V \Lambda $, i.e.~on the interior vertices. 
For the boundary vertices $v\in\bdry V\Lambda $ we set $\tilde S_v:=\{0\}\times\CC^{\deg_{G_\Lambda}}$,
which corresponds to Dirichlet boundary conditions.
This works since in the interior, the degrees with respect to $G$ and $G_\Lambda$ coincide.

Now we turn to the construction of the potential term.  An
\emph{alloy-type potential} is a stochastic process $\cV\colon
\Omega\times G \to \RR$ of the form $\cV_\omega=\sum_{e\in E}
\omega_e \, u_e$, satisfying the following conditions:

The \emph{coupling constants} $\omega_e, e \in E$, form  a family
of independent and identically distributed, non-trivial bounded random variables.
In the operator under consideration each edge $e$ is associated with a single
site potential $u_e$ which is linearly coupled to $\omega_e$.
For this reason the distribution of $\omega_e$ is called a single site
distribution. We denote it by $\mu$.
The expectation of the product measure
$\PP:=\bigotimes_{e\in E} \mu$ is denoted by $\EE$.
Choose $C_\mu$  such that $\supp \mu\subset [-C_\mu,C_\mu]$.

The family of single site potentials $u_e, e\in E$, is assumed to fulfil a partial covering condition and a summability condition:
\begin{dfn}\label{dfn:covering}
Let $I\subset\RR$ be an interval.
The family of single site potentials $u_e, e\in E$, is said to fulfil a \emph{partial covering condition}
with lower bound $c_-(I)>0$ if there is a family of nonempty subintervals $S_e\subset[0,l_e]$ of length $s_e$, $e\in E$,
and for each finite set of edges $\Lambda$ there is a finite set of edges $\Lambda^u$ such that
\begin{align*}
\sum_{e\in \Lambda^u} u_e(x) &\ge c_-(I)  \sum_{e\in\Lambda} C(\lambda,e) \chi_{S_e}, \quad\text{where} \\
C(\lambda,e) &= \frac{l_e}{s_e}\exp\left( 8 \;l_e  \sqrt{C_\mu\| W \|_{L^\infty(e)}+|\lambda|} \right), \\
W &= \sum_{e\in E} u_e \in L^\infty_{loc}(E),
\end{align*}
holds for all $\lambda\in I$.
\end{dfn}
For the following definition, recall that for a metric graph $\tilde G$ with finite set of edges $\tilde E$ and length function $e \mapsto l_e$
the volume is given by $\vol \tilde G=\sum_{e\in\tilde E} l_e$. In contrast to this, $|\Lambda|$ denotes the number of edges in
$\Lambda \subset E$.
\begin{dfn}
Denote by $\Lambda_e$ the minimal set of edges containing the support of $u_e|_\Lambda$
and by $\partial_{\Lambda_e} V$ the boundary vertices of the induced subgraph $G_{\Lambda_e}$.
Then, the family of single site potentials $u_e, e\in E,$ is called \emph{summable} if there are constants $C_j,j=1,2,3$, such that
\begin{equation} \label{eqn:summable}
\begin{aligned}
 \sum_{e\in\Lambda^u} \sum_{v\in \partial_{\Lambda_e} V}\deg v &\leq C_1 |\Lambda| \quad\text{(finite degree property)}, \\
 \sum_{e\in\Lambda^u} \sqrt{\|u_e\|_\infty} \vol G_{\Lambda_e} &\leq C_2 |\Lambda| \quad\text{($L^2$-boundedness)}, \\
 \sum_{e\in\Lambda^u} |\Lambda_e| &\leq C_3 |\Lambda| \quad\text{(volume growth)}
\end{aligned}
\end{equation}
for each finite set of edges $\Lambda$.
\end{dfn}
In particular, this holds if $u_e$ is supported on $e$, uniformly bounded above on $[0,l_e]$ and away from $0$ on $S_e$
(so that $\Lambda^u=\Lambda$, $\Lambda_e=\{e\}$) and there are uniform bounds on vertex degrees and edge lengths.
But our definition is much more general. For instance, decreasing edge lengths can compensate for potential growth and vice versa.

On a given subset $\Lambda\subset E$, we define a
\emph{random Schr\"odinger operator of alloy-type} $H_\omega^\Lambda =-\Delta_\Lambda +\cV_\omega|_{\Lambda}$
as the sum of a negative Laplacian $-\Delta_\Lambda$ on $L^2(G_\Lambda)$ with selfadjoint boundary conditions
$(S_v)_{v\in V_\Lambda}$
and the restriction $\cV_\omega|_{\Lambda}$ of an alloy type potential $\cV_\omega$ to the subgraph $G_\Lambda$.
If $\Lambda=E$, we write $H_\omega$ for $H_\omega^\Lambda$.

\begin{dfn}
The global \emph{modulus of continuity} of the single site distribution $\mu$ is defined
by
\begin{equation}
\label{definition-s-mu-epsilon}
0\le\varepsilon\mapsto s(\mu,\varepsilon):=\sup\{\mu([\lambda-\varepsilon,\lambda +\varepsilon]) \mid \lambda \in \RR\}.
\end{equation}
\end{dfn}
\begin{rem}\label{r:rcmodulus}
$s$ need not be left-continuous as a function of $\varepsilon$ since $\mu$ can have atoms.
However, $s$ is always right-continuous:
\begin{equation} \label{e:rcmodulus}
s(\mu, \varepsilon) = \lim_{\delta \searrow 0} s(\mu, \varepsilon+\delta).
\end{equation}
In fact, assume there is a sequence $\varepsilon_n\searrow\varepsilon$ such that
$\lim_{n\to\infty}s(\mu,\varepsilon_n)=s(\mu,\varepsilon)+2C,C>0$.
Then there exists a sequence $\lambda_n,n\in\NN,$ such that
$\mu([\lambda_n-\varepsilon_n,\lambda_n+\varepsilon_n]) \ge s(\mu,\varepsilon_n)+C$.
Without loss of generality we can assume this sequence to be convergent, since the support of $\mu$ is compact.
Now, the sequence of characteristic functions $f_n:=\chi_{[\lambda_n-\varepsilon_n,\lambda_n+\varepsilon_n]}$
converges to $\chi_{[\lambda-\varepsilon,\lambda+\varepsilon]}$ pointwise at least for $x\ne \lambda-\varepsilon,\lambda+\varepsilon$.
Again, w.l.o.g.\ we can assume this sequence to be convergent everywhere, with limit $0$ or $1$ at $x= \lambda-\varepsilon,\lambda+\varepsilon$,
so that $\lim f_n=:f\leq \chi_{[\lambda-\varepsilon,\lambda+\varepsilon]}$.
Now,  Lebesgue's dominated convergence theorem applies and shows that  $\int f\,d\mu=\lim_n\int f_n\,d\mu\ge s(\mu,\varepsilon)+C$ by construction,
whereas $\int f\,d\mu\le \int \chi_{[\lambda-\varepsilon,\lambda+\varepsilon]}=\mu([\lambda-\varepsilon,\lambda +\varepsilon])\le  s(\mu,\varepsilon)$.

Note also that due to monotonicity in $\varepsilon$, $s(\mu, \varepsilon)\leq s(\mu, t\varepsilon)$ for every $t\ge1$,
and that due to additivity, $s(\mu, M \varepsilon)\leq Ms(\mu, \varepsilon)$ for $M\in \NN$.
Thus, in the following we are free to absorb constant factors appearing in the argument of $s$ into the overall constants in front of $s$.
\end{rem}

With these definitions we can formulate
\begin{thm} \label{t:Wegner}
\label{t:WE} 
Let $\cV_\omega$ be an alloy-type potential  and $\lambda_0 \in \RR$ an energy.
Then there exists a constant $C_W=C_W(\lambda_0)$ such that for all $\lambda\le \lambda_0$,
all finite sets of edges $\Lambda$ and all $\varepsilon \le 1/2$
\begin{equation}
\label{e:WE}
\EE\{\Tr [ \chi_{[\lambda-\varepsilon,\lambda +\varepsilon]}(H_\omega^\Lambda) ]\}
\le C_W \ s(\mu,\varepsilon) \, |\Lambda| \, .
\end{equation}
\end{thm}

\begin{exm}[Random displacements]
Consider a metric graph $G$ as defined above and $u_e, e\in E,$
a collection of single site potentials. Assume that for each $e \in E$ the support $\supp u_e$
as a subset of the interval $[0, l_e]$ satisfies $a_e :=\inf \supp u_e \ge 0$ and $b_e :=\sup \supp u_e \le l_e$.
  Let $(\Omega, \PP), \Omega = \times_{e \in E} \RR, \PP=\otimes_{e\in E} \mu$
and $(\Omega', \PP')$ be two (independent) probability spaces, and $\xi_e \colon \Omega' \to \RR, e \in E,$
a collection of random variables which satisfy
\[
\forall e \in E :\,\,  0 \le \xi_e + a_e \text{ and } \xi_e + b_e
\le l_e.
\]
Denote by $\xi$ the family of random variables $(\xi_e)_{e \in  E}$.
For any $(\omega, \xi) $
define the alloy type potential with random displacements by
\[
 \cV_{\omega,\xi}(x) =\sum_{e\in E} \omega_e u_e (x-\xi_e).
\]
Here we used the convention $u_e(y)=0$, if $y \not\in [0,l_e]$.
Then for any fixed realisation of the vector $\xi$, the stochastic field
$ \cV_{\omega}(x):= \cV_{\omega,\xi}(x)$  defines an alloy type potential satisfying
the requirements of our theorem.
\end{exm}

\begin{exm}[Alloy type potential on a Penrose tiling graph]
Let $G$ be a Penrose tiling graph considered as an embedded metric
graph in $\RR^2$. All its edges have length equal to one. Note that this graph is not periodic.
Assume that the single site potential $u_e$ is bounded, that its support is
contained in the edge $e$, and that there is a subinterval $(a,b)\subset (0,1)$ such that 
$u_e \ge \chi_{(a,b)}$. Then the
potential $ \cV_{\omega}(x):=  \sum_{e\in E} \omega_e u_e (x)$
satisfies the requirements of Theorem~\ref{t:Wegner}.

Similar examples can be given with aperiodic graphs with several different edge-lengths.
\end{exm}

\section{Applications}

There are two important applications of Wegner estimates. Firstly: In combination with
a theorem which establishes the existence and the self-averaging nature of the
IDS it can be used to bound the modulus of continuity of the latter quantity.
Secondly: the Wegner estimate serves as an ingredient in the proof of localisation
via the so called \emph{multiscale analysis}. This method is an induction argument over
increasing length scales and the induction step is based on the Wegner estimate.

\subsection{Modulus of continuity of the IDS}
We describe now how Theorem~\ref{t:Wegner} leads to a continuity
estimate for the IDS for  Schr\"odinger operators on metric graphs
with a $\ZZ^\nu$-structure.

Using the standard embedding of $\ZZ^\nu$ in $\RR^\nu$ we can define:
\begin{dfn}[Metric graphs with $\ZZ^\nu$-structure]
\label{d:GitterGraph}
 A metric graph with $\ZZ^\nu$-structure has $V=\ZZ^\nu$ as vertex
 set and all line segments of length one connecting two vertices
 as edges. The orientation of edges is given by the direction
 of the increasing coordinate.
\end{dfn}
This gives a regular metric graph with degree $2\nu$ and edge lengths $1$.

\begin{dfn}[Alloy type models with $\ZZ^\nu$-structure]
\label{d:GitterModell} An alloy type Schr\"odin\-ger operator on a
metric graph with $\ZZ^\nu$-structure is a random Schr\"odinger
operator of alloy-type $H_\omega, \omega \in \Omega,$ on a metric graph with $\ZZ^\nu$-structure 
(in the sense of Definition~\ref{d:GitterGraph}) such that the boundary
subspace $S_v$ is independent of the vertex $v$ and 
such that for any pair of edges $e, \tilde e$ such that $\tilde{e}=e+k$ 
for some  $k\in\ZZ^\nu$ (in the sense of the embedding $\ZZ^\nu\subset\RR^\nu$), 
we have $u_{\tilde{e}}(\cdot)= u_e(\cdot-k)$.
\end{dfn}
For $l\in\NN$ we define $\Lambda_l$ as the set of edges contained (via the embedding) in $(0,l)^\nu$
which in turn defines the associated subgraph; and we introduce the operator $H_\omega^l$ on the subgraph as the sum of
$-\Delta_{G_{\Lambda_l}}$ and $\cV_\omega|_{G_{\Lambda_l}}$, where $\Delta_{G_{\Lambda_l}}$ is the restriction to $G_\Lambda$ with Dirichlet boundary conditions.

For any $l \in\NN$ and $\omega \in \Omega$ the spectrum of the finite volume operator
$H_\omega^l$ is real, lower bounded and discrete. Thus one may enumerate its eigenvalues
in ascending order and counting multiplicities by
\[
\lambda_1(H_\omega^l) \le
\lambda_2(H_\omega^l) \le \lambda_3(H_\omega^l) \le \dots
\]
For each $\lambda \in \RR$ and $l \in \NN$ the volume-scaled eigenvalue counting function
\[
N_\omega^l(\lambda) := \frac{1}{l^{\nu}} \sharp \{n \in \NN \mid \lambda_n(H_\omega^l) \le \lambda\}
= \frac{1}{l^{\nu}} \Tr \left [ \chi_{(-\infty, \lambda]}(H_\omega^l) \right ]
\]
defines a distribution function. We also define the non-decreasing function
\[
 N(\lambda):=  \sup_{l\in\NN}\EE \{N_\omega^l(\lambda)\}.
\]
Dirichlet-Neumann bracketing arguments show that actually 
\[
\sup_{l\in\NN}\EE \{N_\omega^l(\lambda)\}=\limsup_{l\to\infty}\EE \{N_\omega^l(\lambda)\}.
\]
The following theorem is taken from \cite{HelmV}. It is the analogue of
well known facts for operators on $L^2(\RR^\nu)$ and $\ell^2(\ZZ^\nu)$.

\begin{thm}
\label{t:IDS}
Let $H_\omega, \omega \in \Omega,$ be an alloy type Schr\"odinger operator on a metric graph with
$\ZZ^\nu$-structure as in Definition~\ref{d:GitterModell}.
Then there exists a subset $\Omega' \subset \Omega$ of measure one such that for all $\omega \in \Omega'$ and for all
$\lambda\in \RR$ where $N$ is continuous the convergence
\begin{equation}
\label{e:convergence}
\lim_{l \to \infty} N_\omega^l(\lambda) = N(\lambda)
\end{equation}
holds.
\end{thm}
\begin{cor}
Under the hypotheses of Theorem~\ref{t:IDS},
the IDS obeys for all $0 \le \varepsilon \le 1/2$ the estimate
\[
 N(\lambda+ \varepsilon) - N(\lambda- \varepsilon) \le C_W (\lambda) s(\mu, \varepsilon).
\]
\end{cor}
Here $C_W (\lambda)$ and $ s(\mu,\varepsilon)$ have the same meaning
as in Theorem~\ref{t:Wegner}. In particular, the IDS is continuous
if $s(\mu,0)=0$.
\begin{proof}
The energies where $N$ is continuous are dense in $\RR$.
Assume first that  $\lambda_2 > \lambda_1$ are two such points.
Then
\[
 N(\lambda_2) -N(\lambda_1) = \lim_{l \to\infty} \left(N_\omega^l(\lambda_2) -N_\omega^l(\lambda_1)\right)
\]
for all $\omega\in\Omega'$. Note that the left hand side is non-random and that the terms in the difference on the right hand side
are both uniformly bounded in $\omega$ and $l$.
Thus
\begin{equation}\label{e:N-continuity}
 N(\lambda_2) -N(\lambda_1) = \lim_{l \to\infty} \EE \left\{N_\omega^l(\lambda_2) -N_\omega^l(\lambda_1)\right\}
\le C_W(\lambda_2)\, s\left(\mu,\frac{\lambda_2-\lambda_1}2\right).
\end{equation}
Now let $\lambda \in\RR$ be arbitrary.
By density, there exist sequences $\lambda_{2,n},\,  \lambda_{1,n}, \,  n\in\NN,$
consisting of points of continuity of $N$ such that $\lambda_{2,n}\searrow \lambda+\varepsilon,  \lambda_{1,n} \nearrow \lambda-\varepsilon$
and $\lambda_{2,n}\le \lambda + 1/2,\lambda_{1,n}\ge \lambda - 1/2 $ for all $n \in\NN$. Consequently
\[
 \lim_{n\to \infty} \big(N(\lambda_{2,n}) -N(\lambda_{1,n})\big)
\le \lim_{n\to \infty} C_W(\lambda)\, s\left(\mu,\frac{\lambda_{2,n}-\lambda_{1,n}}2\right)
= C_W(\lambda) s(\mu, \varepsilon)
\]
by property \eqref{e:rcmodulus}. Thus \eqref{e:N-continuity} holds for all $\lambda_1,\lambda_2 \in\RR$.
\end{proof}

\subsection{Localisation via multiscale analysis}
We discuss now the application of Wegner estimates in the proof of localisation based on the multiscale analysis.
This method of proof is applicable to a wide range of random operators.
They have to satisfy certain conditions which can be divided into two groups. The first group of conditions
is rather abstract and can be usually verified for the model at hand without too much effort.
They are commonly stated as follows, see for instance \cite{Stollmann-01}:
\begin{itemize}
\item The underlying space has a $\RR^\nu$ or $\ZZ^\nu$ structure. (Otherwise it may be impossible to define scales.)
 \item The restrictions of the Hamiltonian $H_\omega$ to two finite-volume subsets $\Lambda$ and
$\Lambda'$  are independent random variables, provided $\Lambda$ and $\Lambda'$ are sufficiently far apart.
\item The operator obeys a Weyl type bound, i.e.~for each interval $I\in \RR$ there exists a constant $C$
such that
\[
 \Tr[\chi_{I}(H_\omega^\Lambda)]\le C |\Lambda|.
\]
\item A \emph{geometric resolvent inequality} holds
which relates resolvents of the operator restricted to different finite-volume subsets.
\end{itemize}
For alloy-type random Schr\"odinger operators on $\ZZ^\nu$-metric graphs all these
conditions hold, as can be inferred from \cite{GruberLV,ExnerHS}.
\medskip

There are two more conditions which are more intricate since they depend on the way how randomness enters the model.
% For alloy-type models the only random part of the operator is the potential, and randomness enters
% the potential via a family of i.i.d.~random variables. For such models
% it is crucial whether the random variables have continuous distribution,
% whether they influence the potential in a monoton way,
% and whether  enter in the operator in a linear way.
Here is how these two conditions typically work together to yield localisation:

Let two open energy intervals $I,I_0\subset \RR$ such that $\bar I\subset I_0$ be given.
There exists a length scale $L_0\in\NN$ such that if for any $L_1\in \NN, L_1 \ge L_0$
the conditions
\begin{enumerate}[(H1)]
 \item an \emph{initial length scale estimate} holds for the interval $I$ on scale $L_1$ and
 \item a \emph{weak Wegner estimate} holds on scale $L\in \NN$ for all $\lambda\in I_0$ and $L\ge L_1$
\end{enumerate}
are satisfied, then localisation holds in the energy interval  $I$. As mentioned before, this
means in particular that there is  no continuous spectrum in $I$ and that
eigenfunctions associated to eigenvalues in $I$ decay exponentially, almost surely.

In the above formulation we have used the term \emph{weak Wegner estimate} to
distinguish it from the type of estimate we have established in Theorem~\ref{t:Wegner}.
In many derivations of a Wegner estimate it was assumed that $\mu$ is absolutely continuous
with bounded density, or at least H\"older continuous. It may be asked whether for a localisation proof
it is relevant at all to establish bounds like
\eqref{e:WE} where $s(\mu, \cdot)$ is some unspecified modulus of continuity.
In the following we want to highlight that such estimates have indeed interesting
applications to the multiscale analysis.

We will say that a weak Wegner estimate on the scale $L$  and at the
energy $\lambda$ (with parameters $\beta, q>0$) holds if the
following is true  (cf.~e.g.~\cite{DreifusK-89,Stollmann-01}):
\begin{equation}
\label{e:weakWegner} \PP\left\{\omega \mid \dist
\left(\sigma(H_\omega^L),\lambda\right)\le e^{-L^{\beta}}\right\}
\le \frac{1}{L^q}.
\end{equation}

We call a measure $\mu$ \emph{log-H\"older continuous} with parameter $\alpha>0$,
if for some $c_\mu\in \RR$ and all $0<\varepsilon<1$ we have
\[
s(\mu, \varepsilon) \le c_\mu \frac{1}{|\log \varepsilon|^{\alpha}}.
\]
The following Lemma summarises our observation (which is certainly known among the experts, although we did not find it written up anywhere). 
It holds for random operators
on metric graphs with $\ZZ^\nu$-structure, on the Euclidean space $\RR^\nu$ and on the lattice $\ZZ^\nu$.

\begin{lem}
Let $H_\omega, \omega \in \Omega,$ be a random Schr\"odinger operator and $I_0\subset \RR$ a bounded interval.
Assume that  there exist constants $C_W, L_0$ such that  for all $\varepsilon>0$ and all $\NN\ni L\ge L_0$
\begin{equation*}
\EE\{\Tr [ \chi_{[\lambda-\varepsilon,\lambda +\varepsilon]}(H_\omega^L) ]\}
\le C_W \ s(\mu,\varepsilon) \, L^\nu \, .
\end{equation*}
Assume that the measure $\mu$ is log-H\"older continuous with parameter $\alpha>(q+\nu)/\beta>0$.
Then a weak Wegner estimate with parameters $\beta, q>0$
holds true at all energies $\lambda\in I_0$ and on all scales $\NN \ni L\ge L_1$,
where $L_1:= \max \left(L_0,(C_W c_\mu)^{1/\delta}\right)$ and $\delta:= \alpha\beta -q -\nu>0$.
\end{lem}

\begin{proof}
\begin{align*}
\PP\left\{\omega \mid \dist \left(\sigma(H_\omega^L),\lambda\right)\le e^{-L^{\beta}}\right\}
&\le
\EE\left\{\Tr \left[ \chi_{(\lambda-e^{-L^{\beta}},\lambda+e^{-L^{\beta}})}(H_\omega^L) \right] \right\}
\\
&\le C_W \, L^\nu \, s(\mu, e^{-L^{\beta}})
\\
&\le C_W \, L^\nu \,  c_\mu \frac{1}{|\log e^{-L^{\beta}}|^{\alpha}}
\\
&= C_W \, L^{\nu-\alpha\beta}
\\
&\le L^{-q} \quad \text{ for all } L \ge (C_W c_\mu)^{1/\delta}
\qedhere \end{align*}
\end{proof}

Besides applying to random Schr\"odinger operators on metric graphs by using our Theorem~\ref{t:Wegner},
this Lemma in particular implies that the Wegner estimates of \cite{HundertmarkKNSV-06,CombesHK}
can be used to derive spectral localisation for Anderson and alloy-type operators
on $\ell^2(\ZZ^\nu)$ and $L^2(\RR^\nu)$  whose single site distributions are log-H\"older continuous.
\bigskip

So far our discussion has been quite abstract. Now we
want to give a specific example where hypotheses  (H1) and (H2) are
fulfilled for random operators on metric graphs. This is the case for the model considered in  \cite{ExnerHS} and
\cite{HelmDiss}. There Anderson and strong dynamical localisation  is proved for
an alloy type model with $\ZZ^\nu$-structure. In these
works the single site potentials are characteristic functions
of the edges and the measure $\mu$ is H\"older continuous with support on
an interval $[q_-,q_+],\, 0\le q_-<q_+$. Under these assumptions, in \cite{ExnerHS} a Wegner estimate 
is proved. 
\cite{ExnerHS} provides also an initial length scale estimate. This one, however, relies on a second assumption on the measure which is
sometimes dubbed disorder assumption. This type of technical
assumption is often used if no result on Lifshitz tails is available. The
initial scale estimate derived in \cite{ExnerHS} is the following:
\begin{thm}\label{t-initial}
Assume that the support of $\mu$ is $[q_-,q_+],\, 0\le q_-<q_+$, and that
there exists $\tau>\frac{\nu}{2}$ such that $\mu([q_-,q_-+h])\leq
h^\tau$ for $h$ small.

Then for each $\xi\in (0,2\tau\!-\!\nu)$ there exist
$\beta=\beta(\tau,\xi)\in (0,2)$ and $l^{\ast}=l^{\ast}(\tau,\xi)$
such that
    \begin{equation*}
     \mathbb{P} \big\{\omega \mid \text{dist}\big(\sigma(H_\omega^l), q_-\big) \leq l^{\beta
     -2}\big\}\,\leq\, l^{-\xi}
     \end{equation*}
holds for all $l\geq l^{\ast}$.
\end{thm}

It is possible to combine the results of \cite{ExnerHS,HelmDiss} and the Wegner estimates proved in this paper
to conclude localisation for a wide class of alloy type models on metric graphs with $\ZZ^\nu$-structure:
Assume that $\mu$ is as in Theorem \ref{t-initial}, that it is log-H\"older continuous 
with exponent $\alpha > 4\nu$, and that $u_e \ge \chi_e$. 
In this case we obtain a weak Wegner estimate (for $\nu<q<\frac\alpha2-\nu$) and an initial length scale estimate which 
allow to conclude by the same line of argument as in \cite{HelmDiss} that there is a neighbourhood
of $\inf \sigma (H_\omega)$  where the spectrum is purely localised.

\section{Proofs}

From \cite{GruberLV,GruberV} we will need a lemma on the spectral shift function (SSF for short)
for restrictions to finite subgraphs, whose proof we spell out in more detail below.
For a definition of and more background on the
SSF, see for instance \cite{Kostrykin-99}.
\begin{lem}\label{Potential-Graph}
Let ${G}$ be a finite or infinite metric graph,
${\Lambda}$ a finite subset of its edges,
$-\Delta$ a selfadjoint realisation of the Laplacian on $L^2({G})$ and
$W_1,W_2$ two potentials acting as bounded operators on $L^2({G})$
such that $\supp (W_2-W_1) \subset  G_{\Lambda}$.
Set $H_j=-\Delta +W_j,j=1,2,$ and assume that the SSF $\xi_{H_1,H_2} $ is well defined.
Denote the restriction of $H_j$ to $G_{\Lambda}$ with Dirichlet conditions by $h_j, j=1,2$. Then
\[
\left|\xi_{H_1,H_2} (\lambda) \right|
\le \sum_{v \in \bdry V \Lambda} \deg_G(v) +\left|\xi_{h_1,h_2} (\lambda) \right|.
\]
\end{lem}
\begin{proof}
The basic idea is to decouple the interior of ${G}_{\Lambda}$ from the exterior by choosing appropriate
boundary conditions on $\bdry V\Lambda$.
So, let $H_j^D$ be given by the same differential expression as $H_j$,
but with the domain specified through the following graph-local boundary conditions:
for vertices in $v\in \bdry V\Lambda$ we choose Dirichlet conditions, i.e.\ $S_v:=\{0\}\times\CC^{\deg_{G}(v)}$; for the other vertices of
$G$ we choose the same boundary value subspaces $S_v$ as for $H_j$.
Then
\[ \left| \xi_{H_j, H_j^D}(\lambda)\right| \leq \sum_{v \in \bdry V\Lambda} \deg_G(v) \]
according to Corollary 11 and Lemma 13 of \cite{GruberLV}, since this is the maximal rank of the perturbation induced
by changing boundary conditions on $\bdry V\Lambda$.

Now, by construction the subgraphs $G_\Lambda$ and $G_{E\setminus\Lambda}$ of $G$ have disjoint complementing edge sets $\Lambda$ and $E\setminus\Lambda$,
their sets of interior vertices $\inter V\Lambda$ and $\inter V{E\setminus\Lambda}=\exter V\Lambda$ are disjoint,
and the intersection of their vertex sets is precisely $\bdry V\Lambda=\bdry V{E\setminus\Lambda}$,
which are the only vertices which connect these subgraphs to each other.
Since we chose Dirichlet conditions for $H_j^D$ on these connecting vertices and since Dirichlet conditions decouple,
the $H_j^D$ decompose into a direct sum of interior and exterior parts, i.e.\ operators on $G_\Lambda$ and $G_{E\setminus\Lambda}$.
The latter coincide by assumption,
the former are given by $h_j$. This proves the assertion.
\end{proof}

From \cite{GruberLV,GruberV} we also need a lemma on the spectral shift for potential perturbations:
\begin{lem}\label{Potential-Kante}
Let $-\Delta$ be a selfadjoint realisation of the Laplacian on an arbitrary finite metric graph $ G$
and let  $W_1,W_2$ be bounded potentials on $L^2( G)$.
Set $H_j=-\Delta +W_j,j=1,2$.
Then
\[
\left|\xi_{H_1,H_2} (\lambda) \right| \leq \left(\sqrt{\|W_1\|}+\sqrt{\|W_2\|}\right) \frac{\vol  G}\pi  +5|E( G)|
\]
where $\vol  G= \sum_{e\in E( G)} l_e$ is the one-dimensional volume of $ G$.
\end{lem}

We need  one more preparatory lemma before giving the proof of Theorem~\ref{t:Wegner}.

Let $\rho$ be a smooth, monotone switch function
$\rho:=\rho_{\lambda,\varepsilon}\colon \RR \to [-1,0]$. \label{p-rho} By a switch function we mean that
for a positive $\varepsilon \le 1/2$, $\rho$  has the following properties: $ \rho\equiv -1$ on
$(-\infty,\lambda-\varepsilon]$, $\rho\equiv 0$ on $[\lambda+\varepsilon,\infty)$ and $\|\rho'\|_\infty \le
1/\varepsilon$. Set $I=[\lambda_0-1,\lambda_0+1]$.
\begin{lem}\label{Traeger-Reduktion}
Under the assumptions of Theorem~\ref{t:Wegner} there is a constant $C_{uc}>0$ depending only on $I$ and $u$
such that
\[
\Tr  [\rho (H_\omega^\Lambda+ C_{uc}\varepsilon)] \le \Tr\Big[\rho (H_\omega^\Lambda+ \varepsilon\sum_{e\in \Lambda^u} u_e) \Big] .
\]
\end{lem}
\begin{proof}
We know from \cite{KirschV-02b,HelmV} (keeping track of good Sobolev constants and lengths)
that for an eigenfunction $\psi$ to the eigenvalue $\lambda$
\begin{align*}
\int_{S_e} |\psi_e|^2 &\geq C(\lambda,e)^{-1} \int_0^{l_e} |\psi_e|^2 ,\quad\text{where} \\
 C(\lambda,e) &= \frac{l_e}{s_e}\exp\left( 8 \;l_e  \sqrt{C_\mu\| W \|_{L^\infty(e)}+|\lambda|} \right).
\end{align*}
Therefore,
\begin{align*}
\sum_{e\in\Lambda} C(\lambda,e) \int_{S_e} |\psi_e|^2 &\geq \sum_{e\in\Lambda} \int_0^{l_e} |\psi_e|^2 =1
\end{align*}
by normalisation. Denote by $\lambda_n^\Lambda(\omega)$ the $n$-th eigenvalue of $H_\omega^\Lambda$
and by $\psi_n$ the associated, normalised eigenfunction.
Since $\sum_{e \in \Lambda^u} \frac{\partial \lambda_n^\Lambda(\omega)}{\partial \omega_e} = \sum_{e \in \Lambda^u} (\psi_n , u_e\psi_n)$,
our partial covering condition from Definition~\ref{dfn:covering} and the explicit form of the estimate above finally lead to
\begin{equation}
\label{Eder}
\sum_{e \in \Lambda^u} \frac{\partial \lambda_n^\Lambda(\omega)}{\partial \omega_e}
\ge C_1(I) >0
\end{equation}
for all eigenvalues $\lambda_n^\Lambda$ of $H_\omega^\Lambda$
inside a bounded energy interval $I$. The bound $C_1(I)$ does not
depend on the set of edges  $\Lambda \subset E$ and on the eigenvalue index
$n\in \NN$.
Set $ \cV_{\omega,t}:= \cV_\omega + t \sum_{e \in \Lambda_u}  u_e $, and denote by
$\lambda_n^\Lambda(\omega,t)$ the eigenvalues of the corresponding Schr\"odinger operator restricted to the finite graph $G_\Lambda$.
Then
\begin{equation*}
\frac{\partial \lambda_n^\Lambda(\omega,t)}{\partial t}
= \sum_{e \in \Lambda^u} \frac{\partial \lambda_n^\Lambda(\omega)}{\partial \omega_e}
\ge C_1(I).
\end{equation*}
Thus
\[
\lambda_n^\Lambda(\omega,\varepsilon) -\lambda_n^\Lambda(\omega,0)
=\int_0^\varepsilon \frac{\partial \lambda_n^\Lambda(\omega,t)}{\partial t}  dt
\ge C_{uc} \varepsilon.
\]
By the isotonicity of $\rho$ it follows that
\[
 \rho(\lambda_n^\Lambda(\omega,\varepsilon)) \ge \rho( \lambda_n^\Lambda(\omega,0)  + C_{uc} \varepsilon).
\qedhere \]

\end{proof}

\begin{proof}[Proof of Theorem~\ref{t:Wegner}]
According to our choice of switch function $\rho$
\[
\chi_{[\lambda-\varepsilon,\lambda +\varepsilon]} (x) \le \rho(x+2\varepsilon) -\rho(x-2\varepsilon).
\]
By Lemma~\ref{Traeger-Reduktion} we conclude (putting $\varepsilon'=\varepsilon/C_{uc}$)
\[
\Tr  [\rho (H_\omega^\Lambda+ \varepsilon)] \le \Tr\Big[\rho (H_\omega^\Lambda+ \varepsilon'\sum_e u_e) \Big] .
\]
Let $\Lambda^u$ be as above.
$\Lambda^u$ contains $L:=|\Lambda^u|$ edges. We
enumerate the edges in $\Lambda^u$ by $ e\colon \{1, \dots, L\} \to \Lambda^u $, $n\mapsto e(n)$, and set
\[
W_0 \equiv 0, \quad  W_n =\sum_{m=1}^{n} u_{e(m)}, \qquad n=1,2,\dots, L.
\]
Thus
\begin{align}
\chi_{[\lambda-\varepsilon,\lambda +\varepsilon]} (H_\omega^\Lambda)
& \le       \nonumber
 \rho(H_\omega^\Lambda+2\varepsilon)-\rho(H_\omega^\Lambda-2\varepsilon)
\\      \label{e:project}
& \le
\rho(H_\omega^\Lambda-2\varepsilon+4\varepsilon' W_{L})-\rho(H_\omega^\Lambda-2\varepsilon)
\\      \nonumber
& =
\sum_{n=1}^{L} \rho(H_\omega^\Lambda-2\varepsilon+4\varepsilon'  W_{n})-
\rho(H_\omega^\Lambda-2\varepsilon+4\varepsilon'  W_{n-1}).
\end{align}
We fix $n \in \{1, \dots, L\}$, define
\[
\omega^\perp := \{\omega_e^\perp\}_{e \in \Lambda^u}, \qquad
\omega_e^\perp :=\begin{cases} 0 \quad &\text{if }e=e(n), \\
\omega_e \quad &\text{if } e\neq e(n), \end{cases}
\]
and set
\[
\phi_{n}(\eta) = \Tr\bigl[\rho(H_{\omega^\perp}^\Lambda-2\varepsilon +4\varepsilon' W_{n-1}+\eta
u_{e(n)})\bigr], \quad \eta\in\RR.
\]
The function $\phi_{n}$ is continuously differentiable, monotone increasing and bounded.
By definition of $\phi_{n}$,
\[
\Tr [ \rho(H_\omega^\Lambda-2\varepsilon+4\varepsilon'  W_n)
-\rho(H_\omega^\Lambda-2\varepsilon+4\varepsilon'  W_{n-1})]
=
\phi_{n}(\omega_{e(n)}+4\varepsilon')-\phi_{n}(\omega_{e(n)})
\]
since $\phi_{n}(\eta) = \Tr\bigl[\rho(H_{\omega}^\Lambda-2\varepsilon +4\varepsilon' W_{n-1}+(\eta
-\omega_{e(n)})u_{e(n)})\bigr]$,
so that
\begin{multline*}
\EE_{\omega_{e(n)}} \{ \Tr [ \rho(H_\omega^\Lambda-2\varepsilon+4\varepsilon'  W_n)
-\rho(H_\omega^\Lambda-2\varepsilon+4\varepsilon'  W_{n-1})] \}
= \\
\int [\phi_{n}(\omega_{e(n)}+4\varepsilon')-\phi_{n}(\omega_{e(n)})]\, d\mu(\omega_{e(n)})
\end{multline*}
where $\EE_{\omega_{e(n)}}$ denotes the expectation with respect to the random variable $\omega_{e(n)}$ only.
Let $\supp(\mu)\subset(a,b)$.
Using Lemma 6 in \cite{HundertmarkKNSV-06}
 we have
\begin{align*}
\int [\phi_{n}(\omega_{e(n)}+4\varepsilon')-\phi_{n}(\omega_{e(n)})]\, d\mu(\omega_{e(n)})
&\le
s(\mu,2\varepsilon') [\phi_{n}(b+4\varepsilon')-\phi_{n}(a)].
\end{align*}
Denote by $\xi_{\Lambda,n}$ the SSF associated to the pair of operators $H_n(a),H_n(b+2/C_{uc})$ on $L^2(\Lambda)$
where $H_n(\eta)$ is given by $H_n(\eta):=H_\omega^\Lambda-2\varepsilon+4\varepsilon'  W_{n-1}+ (\eta-\omega_{e(n)}) u_{e(n)} $.
Then by the Krein trace identity and the normalisation of $\rho$
\begin{align*}
\phi_{n}(b+4\varepsilon')-\phi_{n}(a)= \int_a^{b+4\varepsilon'} \rho' \ \xi_{\Lambda,n} \ d\lambda
\le \|\xi_{\Lambda,n}\|_\infty.
\end{align*}
Let $\Lambda_e$, $u_e$, $\bdry V{\Lambda_e}$ and $G_{\Lambda_e} $ be as in the definition of summable potentials.
By $\xi_{\Lambda_{e(n)},n}$ we denote the SSF associated to the pair $H_n(a),H_n(b+2/C_{uc})$,
but now considered as operators on $L^2(\Lambda_{e(n)})$.
Apply Lemma~\ref{Potential-Graph} to obtain:
\begin{align*}
\|\xi_{\Lambda,n}\|_\infty
\le \sum_{v \in \bdry V{\Lambda_{e(n)}}} \deg_G(v) +\|\xi_{\Lambda_{e(n)},n}\|_\infty
\end{align*}
Now apply Lemma~\ref{Potential-Kante} successively $L$ times to obtain
\begin{multline*}
\EE \{\Tr [\chi_{[\lambda-\varepsilon,\lambda +\varepsilon]} (H_\omega^\Lambda) ]\}
\\
\begin{aligned}
& \le s(\mu,2\varepsilon') \;  \sum_{n=1}^{L} \left( \sum_{v \in \bdry V{\Lambda_{e(n)}}} \deg_G(v)
     + \sqrt{\|u_{e(n)}\|_\infty}\; \frac{\vol G_{\Lambda_{e(n)}}}\pi +  5 | \Lambda_{e(n)}| \right) \\
&\le s(\mu,2\varepsilon') \; (C_1+C_2/\pi+5C_3) |\Lambda|
\end{aligned}
\end{multline*}
by the summability condition on the family of single site potentials.
\end{proof}

  \bibliographystyle{amsalpha}	
 \bibliography{ghv}

\end{document}